\overfullrule=0pt
\centerline {\bf A note on nonlocal discrete problems involving sign-changing Kirchhoff functions}\par
\bigskip
\bigskip
\centerline {BIAGIO RICCERI}\par
\bigskip
\bigskip
{\bf Abstract.} In this note, we establish a multiplicity theorem for a nonlocal discrete problem of the type
$$\cases{-\left(a\sum_{m=1}^{n+1}|x_m-x_{m-1}|^2+b\right)(x_{k+1}-2x_k+x_{k-1})=h_k(x_k)\hskip 10pt k=1,...,n, \cr & \cr
x_0=x_{n+1}=0\cr}$$
assuming $a>0$ and (for the first time) $b<0$.
\par
\bigskip
\bigskip
\bigskip
\bigskip
{\bf 1. Introduction}\par
\bigskip
If $\Omega$ is a bounded domain of ${\bf R}^n$ and $K:[0,+\infty[\to {\bf R}$, $\varphi:\Omega\times {\bf R}\to {\bf R}$ are two given functions, the nonlocal problem
$$\cases {-K(\int_{\Omega}|\nabla u(x)|^2)\Delta u=
\varphi(x,u)
 & in
$\Omega$\cr & \cr 
u=0 & on
$\partial \Omega$\cr}$$
 is certainly among the most studied ones in today's nonlinear analysis (we refer to [7] for an introduction to the subject).\par
\smallskip
In checking the relevant literature, one can realize 
that, in the majority of the papers, one assumes $K(t)=at+b$ with $a>0$ and $b\geq 0$ and, in any case, that
the Kirchhoff function $K$ is assumed to be, in particular, continuous and non-negative in $[0,+\infty[$.\par
\smallskip
However, it is natural to ask
what happens when at least one of these properties fails.\par
\smallskip
The case where $K$ can be discontinuous in $[0,+\infty[$ has been considered for the first
time in [4], for $n=1$, and then in [5] for the general case (see also [1], [2], [3]). In these papers, however, $K$ is non-negative.\par
\smallskip
 The papers dealing with
a sign-changing function $K$ are more numerous, but in each of them it is assumed that $K(t)=at+b$ with $a<0$ and $b\geq 0$. The first of these
papers was [8].\par
\smallskip
In the present very short note, we are interested in the discrete counterpart of the above problem. That is to say, given $n$ continuous functions
$f_k:{\bf R}\to {\bf R}$ ($k=1,...,n)$, we deal with the problem
$$\cases{-K\left(\sum_{h=1}^{n+1}|x_h-x_{h-1}|^2\right)(x_{k+1}-2x_k+x_{k-1})=f_k(x_k)\hskip 10pt k=1,...,n, \cr & \cr
x_0=x_{n+1}=0\ . \cr}$$
Also for this discrete problem, we can repeat what we said before, even in a stronger way: it seems that in each paper on the subject, the
function $K$ is continuous and non-negative in $[0,+\infty[$.\par
\smallskip
Our aim is to establish a multiplicity result for this problem where (for the first time) the Kirchhoff function $K$ changes sign.\par
\bigskip
{\bf 2. Results}\par
\bigskip
Before stating our result, we recall the following two theorems which will be key tools used in our proof. 
\bigskip
THEOREM 2.A.([6]) - {\it Let $X$ be a topological space, let $Y$ be a convex set in a topological vector space and let $h:X\times Y\to {\bf R}$
 be lower semicontinuous and inf-compact in $X$, and continuous and quasi-concave in $Y$. Also, assume that
$$\sup_Y\inf_Xh<\inf_X\sup_Yh\ .$$
Moreover, let $\varphi:X\to {\bf R}$ be a lower semicontinuous function such that
$$\sup_X\varphi-\inf_X\varphi<\inf_X\sup_Yh-\sup_Y\inf_Xh.$$
Then, for each convex set $S\subseteq Y$, dense in $Y$, there exists $\tilde y\in S$ such that the function $h(\cdot,\tilde y)+\varphi(\cdot)$ has at least two global minima.}\par
\bigskip
THEOREM 2.B.([6]) - {\it Let $X$ be a topological space, let $H$ be a real Hilbert, let $Y$ be a closed ball in $H$ centered at $0$,
and let $Q:X\to {\bf R}$, $\psi:X\to H$. Assume that the functional $x\to Q(x)-\langle\psi(x),y\rangle$ is lower semicontinuous
for each $y\in Y$, while the functional $x\to Q(x)-\langle\psi(x),y_0\rangle$ is inf-compact for some $y_0\in Y$. Moreover,
assume that, for each $x\in X$, there exists $u\in X$ such that
$$Q(x)=Q(u)$$
and
$$\psi(x)=-\psi(u)\ .$$
Finally, assume that there is no global minimum of $Q$ at which $\psi$ vanishes.\par
Then, we have
$$\sup_{y\in Y}\inf_{x\in X}(Q(x)-\langle\psi(x),y\rangle)<\inf_{x\in X}\sup_{y\in Y}(Q(x)-\langle\psi(x),y\rangle)\ .$$}\par
\medskip
Our main result is as follows:\par
\medskip
THEOREM 2.1. - {\it Let $K:[0,+\infty[\to {\bf R}$, $f_1,...,f_n:{\bf R}\to {\bf R}$ be $n+1$ continuous functions satisfying
the following conditions:\par
\noindent
$(a)$\hskip 5pt $\inf_{t>0}\int_0^tK(s)ds<0$ and $\liminf_{t\to +\infty}{{\int_0^tK(s)ds}\over {t}}>0$;\par
\noindent
$(b)$\hskip 5pt $\limsup_{|t|\to +\infty}{{\left|\int_0^tf_k(s)ds\right|}\over {t^2}}<+\infty$ for each $k=1,...,n$;\par
\noindent
$(c)$\hskip 5pt for each $k=1,...,n$, the function $t\to \int_0^tf_k(s)ds$ is odd and vanishes only at $0$.\par
Then, for each $r>0$, there exists a number $\delta>0$ with the following property: for every n-uple of continuous functions $g_1,...,g_n:{\bf R}\to {\bf R}$
satisfying
$$\max_{1\leq k\leq n}\left(\sup_{t\in {\bf R}}\int_0^tg_k(s)ds-\inf_{t\in {\bf R}}\int_0^tg_k(s)ds\right)<\delta,$$
there exists $(\tilde\mu_1,...\tilde\mu_n)\in {\bf R}^n$, with $\sum_{k=1}^n|\tilde\mu_k|^2<r^2$, such that the problem
$$\cases{-K\left(\sum_{h=1}^{n+1}|x_h-x_{h-1}|^2\right)(x_{k+1}-2x_k+x_{k-1})=g_k(x_k)+\tilde \mu_kf_k(x_k)\hskip 10pt k=1,...,n, \cr & \cr
x_0=x_{n+1}=0\cr}$$
has at least three solutions.}\par
\smallskip
PROOF. Fix $r>0$. First, we are going to apply Theorem 2.B. In this connection,  take
$$X=\{(x_0,x_1,...,x_n,x_{n+1})\in {\bf R}^{n+2} : x_0=x_{n+1}=0\},$$
with the scalar product
$$\langle x,y\rangle _1=\sum_{k=1}^{n+1}(x_k-x_{k-1})(y_k-y_{k-1}).$$
We denote by $\langle\cdot,\cdot\rangle_2$ the usual scalar product on ${\bf R}^n$. That is
$$\langle x,y\rangle_2=\sum_{k=1}^nx_ky_k.$$
Fix $\gamma>0$ so that
$$\|x\|_2\leq \gamma\|x\|_1 \eqno{(2.1)}$$
for all $x\in X$.
Consider the functions
$Q:X\to {\bf R}$ and $\psi:X\to {\bf R}^n$ defined by
$$Q(x)={{1}\over {2}}\int_0^{\|x\|_1^2}K(s)ds$$
and 
$$\psi(x)=\left(\int_0^{x_1}f_1(s)ds,...,\int_0^{x_n}f_n(s)ds\right)$$
for all $x\in X$. Fix $\mu\in {\bf R}^n$. In view of $(a)$ and $(b)$, there exists
$\eta_1, \eta_2,\eta_3>0$ such that
$$\int_0^tK(s)ds\geq\eta_1t-\eta_2 \eqno{(2.2)}$$
for all $t\geq 0$ 
and
$$\left|\int_0^{t}f_k(s)ds\right|\leq \eta_3(t^2+1) \eqno{(2.3)}$$
for all $t\in {\bf R}$, $k=1,...,n$. Fix $x\in X$. Using $(2.2)$ and the Cauchy-Schwarz inequality, we obtain
$$Q(x)-\langle\psi(x),\mu\rangle_2
\geq {{1}\over {2}}\eta_1\|x\|_1^{2}-{{1}\over {2}}\eta_2-|\langle\psi(x),\mu\rangle_2|\geq 
{{1}\over {2}}\eta_1\|x\|_1^{2}-{{1}\over {2}}\eta_2-\|\mu\|_2\|\psi(x)\|_2. \eqno{(2.4)}$$
On the other hand, in view of $(2.3)$,  for each $k=1,...,n$, we have
$$\left|\int_0^{x_k}f_k(s)ds\right|\leq \eta_3(|x_k|^2+1)$$
and hence
$$\|\psi(x)\|_2\leq \eta_3\sqrt{\sum_{k=1}^n(|x_k|^2+1)^2}\leq \eta_3\left(\sum_{k=1}^n|x_k|^2+n\right). \eqno{(2.5)}
$$
Putting $(2.1), (2.4)$ and $(2.5)$ together, we get
$$Q(x)-\langle\psi(x),\mu\rangle_2\geq {{1}\over {2}}\eta_1\|x\|_1^{2}-\|\mu\|_2\eta_3(\gamma^2\|x\|_1^2+n)-{{1}\over {2}}\eta_2.
\eqno{(2.6)}$$
Now, fix $\sigma>0$ so that
$$\sigma<\min\left\{{{\eta_1}\over {2\eta_3\gamma^2}},r\right\}.$$
Let $Y$ be the closed ball in ${\bf R}^n$ centered at $0$, of radius $\sigma$. If $\mu\in Y$, in view of $(2.6)$, we have
$$\lim_{\|x\|_1\to +\infty}(Q(x)-\langle\psi(x),\mu\rangle_2)=+\infty$$
and so the function $x\to Q(x)-\langle\psi(x),\mu\rangle_2$ is inf-compact.
Further, observe that, by $(c)$, the function $\psi$ vanishes only at $0$, while, by $(a)$, $0$ is not a global minimum of $Q$. Clearly, $Q$ is
even and $\psi$ is odd, in view of $(c)$. In other words, each assumption of Theorem 2.B is satisfied. Consequently, the number
$$\delta:={{1}\over {n}}\left(\inf_{x\in X}\sup_{\mu\in Y}(Q(x)-\langle\psi(x),\mu\rangle_2)-\sup_{\mu\in Y}\inf_{x\in X}(Q(x)-\langle\psi(x),\mu\rangle_2)\right)\eqno{(2.7)}$$
is positive. At this point, we apply Theorem 2.A taking
$$h(x,\mu)=Q(x)-\langle\psi(x),\mu\rangle_2$$
for all $(x,\mu)\in X\times Y$. Fix $n$  continuous functions $g_1,...,g_n:{\bf R}\to {\bf R}$
satisfying
$$\max_{1\leq k\leq n}\left(\sup_{t\in {\bf R}}\int_0^tg_k(s)ds-\inf_{t\in {\bf R}}\int_0^tg_k(s)ds\right)<\delta\eqno{(2.8)}$$
and consider the function $\varphi:X\to {\bf R}$
defined by
$$\varphi(x)=-\sum_{k=1}^n\int_0^{x_k}g_k(s)ds$$
for all $x\in X$. Clearly, in view of $(2.7)$ and $(2.8)$, we have
$$\sup_X\varphi-\inf_X\varphi\leq \sum_{k=1}^n\left(\sup_{t\in {\bf R}}\int_0^tg_k(s)ds-\inf_{t\in {\bf R}}\int_0^tg_k(s)ds\right)$$
$$\leq
n\max_{1\leq k\leq n}\left(\sup_{t\in {\bf R}}\int_0^tg_k(s)ds-\inf_{t\in {\bf R}}\int_0^tg_k(s)ds\right)<\inf_X\sup_Yh-\sup_Y\inf_Xh.
$$
So, each assumption of Theorem 2.A is satisfied. As a consequence, there exists $\tilde \mu\in Y$ such that the function
$$J_{\tilde\mu}(\cdot):=h(\cdot,\tilde \mu)+\varphi(\cdot)$$
has at least two global minima in $X$. 
It is clear that this function $J_{\tilde\mu}$ is $C^1$, with derivative given by
$$J_{\tilde\mu}'(x)(y)=K\left(\sum_{h=1}^{n+1}|x_h-x_{h-1}|^2\right)\langle x,y\rangle_1-\sum_{k=1}^ng_k(x_k)y_k-\sum_{k=1}^n\tilde\mu_kf_k(x_k)y_k$$
for all $x,y\in X$. So, taking into account that
$$\langle x,y\rangle_1=-\sum_{k=1}^n(x_{k+1}-2x_k+x_{k-1})y_k\ ,$$
we have
$$J_{\tilde\mu}'(x)(y)=-K\left(\sum_{h=1}^{n+1}|x_h-x_{h-1}|^2\right)\sum_{k=1}^n(x_{k+1}-2x_k+x_{k-1})y_k-\sum_{k=1}^ng_k(x_k)y_k-\sum_{k=1}^n\tilde\mu_kf_k(x_k)y_k
\eqno{(2.9)}$$
for all $x,y\in X$. Since $J_{\tilde\mu}$ is coercive and has at least two global minima, by a classical theorem of Courant, it possesses at least three critical points which,
by $(2.9)$, are three solutions of the problem.\hfill $\bigtriangleup$
\medskip
Here is a remarkable corollary of Theorem 2.1:\par
\medskip
COROLLARY 2.1. - {\it Let $f_1,...f_n:{\bf R}\to {\bf R}$ be $n$ continuous functions satisfying conditions $(b)$ and $(c)$ of Theorem 2.1.\par
Then, for each $a,r>0$ and $b<0$, there exists a number $\delta>0$ with the following property: for every n-uple of continuous functions $g_1,...,g_n:{\bf R}\to {\bf R}$
satisfying
$$\max_{1\leq k\leq n}\left(\sup_{t\in {\bf R}}\int_0^tg_k(s)ds-\inf_{t\in {\bf R}}\int_0^tg_k(s)ds\right)<\delta,$$
there exists $(\tilde\mu_1,...\tilde\mu_n)\in {\bf R}^n$, with $\sum_{k=1}^n|\tilde\mu_k|^2<r^2$, such that the problem
$$\cases{-\left(a\sum_{h=1}^{n+1}|x_h-x_{h-1}|^2+b\right)(x_{k+1}-2x_k+x_{k-1})=g_k(x_k)+\tilde \mu_kf_k(x_k)\hskip 10pt k=1,...,n, \cr & \cr
x_0=x_{n+1}=0\cr}$$
has at least three solutions.}\par
\smallskip
PROOF. It is enough to observe that the function $K(t)=at+b$ satisfies condition $(a)$ of Theorem 2.1.\hfill $\bigtriangleup$\par
\medskip
REMARK 1.1. - It is important to remark that the technique adopted in the proof Theorem 2.1 cannot be used to treat the non-discrete problem,
keeping condition $(a)$. This is due to the fact that, under condition $(a)$, 
the functional $$u\to \int_{0}^ {\int_{\Omega}|\nabla u(x)|^2dx}K(s)ds$$
 is not weakly lower semicontinuous in $H^1_0(\Omega)$.
\bigskip
\bigskip
{\bf Acknowledgements.} This work has been funded by the European Union - NextGenerationEU Mission 4 - Component 2 - Investment 1.1 under the Italian Ministry of University and Research (MUR) programme "PRIN 2022" - grant number 2022BCFHN2 - Advanced theoretical aspects in PDEs and their applications - CUP: E53D23005650006. The author has also been supported by the Gruppo Nazionale per l'Analisi Matematica, la Probabilit\`a e 
le loro Applicazioni (GNAMPA) of the Istituto Nazionale di Alta Matematica (INdAM) and by the Universit\`a degli Studi di Catania, PIACERI 2024-2026, Linea di intervento 2, Progetto ”PAFA”. \par
\bigskip
\bigskip
\centerline {\bf References}\par
\bigskip
\bigskip
\noindent
[1]\hskip 5pt Y. H. KIM, {\it Existence and uniqueness of solutions to non-local problems of Brezis-Oswald type ´
and its application} Fractal Fract., 8 (2024), 622. https://doi.org/10.3390/fractalfract8110622
\smallskip
\noindent
[2]\hskip 5pt Y. H. KIM, {\it Existence and uniqueness of solution to the p-Laplacian equations involving discontinuous Kirchhoff functions via a global minimum principle of Ricceri}, Minimax Theory Appl., {\bf 10} (2025), 34-42.\par
\smallskip
\noindent
[3]\hskip 5pt I. H. KIM and Y. H. KIM,  {\it Existence, uniqueness, and localization of positive solutions to nonlocal problems of the Kirchhoff type via the global minimum principle of Ricceri}, AIMS Mathematics, 2025, 10(3): 4540-4557.\par
\smallskip
\noindent
[4]\hskip 5pt B. RICCERI, {\it Multiplicity theorems involving functions with non-convex range}, Stud. Univ. Babe\c{s}-Bolyai Math., {\bf 68} (2023), 125-137.\par
\smallskip
\noindent
[5]\hskip 5pt B. RICCERI, {\it Existence, uniqueness, localization and minimization property of positive solutions for non-local problems involving discontinuous Kirchhoff functions}, Adv. Nonlinear Anal., {\bf 13} (2024), Paper No. 20230104. \par
\smallskip
\noindent
[6]\hskip 5pt B. RICCERI, {\it Multiple critical points in closed sets via minimax theorems}, Optimization (2025),\par
\noindent
https://doi.org/10.1080/02331934.2025.2457549.
\smallskip
\noindent
[7]\hskip 5pt P. PUCCI and V. D. R\u ADULESCU, {\it Progress in nonlinear Kirchhoff problems}, Nonlinear Anal., {\bf 186} (2019), 1-5.\par
\smallskip
\noindent
[8]\hskip 5pt G. YIN and J. LIU, {\it Existence and multiplicity of nontrivial solutions for a nonlocal problem},  
Bound. Value Probl., 2015, 2015:26, 1-7.\par
\smallskip
\noindent

\bye